\documentclass{article}

\usepackage{amsmath,amssymb}

\input xy
\xyoption{all}

\def\ad{\hbox{\rm ad}\ }
\def\hxl{{co-split Lie}}
\def\hx{{co-splitting}\ }
\def\mc{\mathbf C}
\def\mr{\mathbf R}
\def\id{\hbox{\rm id}}
\def\beqs{\begin{eqnarray*}}
\def\eeqs{\end{eqnarray*}}
\def\beq{\begin{eqnarray}}
\def\eeq{\end{eqnarray}}

\begin{document}

\title{Introduction to co-split Lie algebras}

\author{Limeng Xia\\
\it\small Faculty of Science, Jiangsu University\\
\it\small Zhenjiang 212013, Jiangsu Province, P.R. China\\
Naihong Hu\footnote{Corresponding author:
nhhu@math.ecnu.edu.cn.}\\
\it\small Department of Mathematics, East China Normal University\\
\it\small Minhang Campus, Dongchuan Road 500, Shanghai 200241, P.R.
China}

%Mathematics Subject Classifications (2000): 17B62

%Keywords: co-split Lie algebra£¬Lie coalgebra, adjoint action

\date{}

\newtheorem{theo}{Theorem}[section]
\newtheorem{defi}{Definition}[section]
\newtheorem{rem}{Remark}[section]
\newtheorem{prop}{Proposition}[section]
\newtheorem{lemm}{Lemma}[section]

\maketitle
\begin{abstract}
In this work, we introduce a new concept which is obtained by
defining a new compatibility condition between Lie algebras and Lie
coalgebras. With this terminology, we describe the interrelation
between the Killing form and the adjoint representation in a new
perspective.
\end{abstract}

\section{Introduction}

During the past decade, there have appeared a number of papers on
the study of Lie bialgebras (see \cite{EK}, \cite{ES} and references
therein, etc). It is well-known that a Lie bialgebra is a vector
space endowed simultaneously with a Lie algebra structure and a Lie
coalgebra structure, together with a certain compatibility
condition, which was suggested by a study of Hamiltonian mechanics
and Poisson Lie groups (\cite{ES}).

In the present work, we consider a new [Lie algebra]-[Lie coalgebra]
structure, say, a \hxl algebra. Using this concept, we can easily
study the Lie algebra structure on the dual space of a semi-simple
Lie algebra from another point of view.

This paper is arranged as follows: At first we recall some concepts
and study the relations between Lie algebras and Lie coalgebras.
Then  we give the definition of a co-split Lie algebra. In section
4, we prove that $sl_{n+1}(\mc)$ is a co-split Lie algebra. Then we
discuss the interrelation of the Killing form and the adjoint
representation of $sl_{n+1}(\mc)$. Finally, the results are proved
to hold for all finite dimensional complex semi-simple Lie algebras.

\section{Basics}

In this section, we mainly recall the definitions of Lie algebras,
Lie coalgebras and Lie bialgebras, and also their relationship. For
more information, one can see \cite{EK}, \cite{ES} and references
therein.

A Lie algebra is a pair $(L,[,])$, where $L$ is a linear space and
$[,]: L\times L\longrightarrow L$ is a bilinear map (in fact, it is
a linear map from $L\otimes L$ to $L$) satisfying:

\quad$\begin{tabular}{cl}
(L1)&[a,b]+[b,a]=0,\\
(L2)&[a,[b,c]]+[b,[c,a]]+[c,[a,b]]=0.
\end{tabular}$

For any spaces $U,V,W$, define maps \beqs \tau:&&U\otimes
V\longrightarrow V\otimes U\\
&&u\otimes v\longmapsto v\otimes u,\\
\xi:&&U\otimes V\otimes W\longrightarrow V\otimes W\otimes
U\\
&&u\otimes v\otimes w\longmapsto v\otimes w\otimes u.\eeqs

A Lie coalgebra is a pair $(L,\delta)$, where $L$ is a linear
space and $\delta:L\longrightarrow L\otimes L$ is a linear map
satisfying:

\quad \begin{tabular}{cl}
(Lc1)&$(1+\tau)\circ\delta=0$,\\
(Lc2)&$(1+\xi+\xi^2)\circ(1\otimes\delta)\circ\delta=0$.
\end{tabular}

A Lie bialgebra is a triple $(L, [,], \delta)$ such that

\quad\begin{tabular}{cl} (Lb1)& $(L,[,])$ is a Lie algebra,\\
(Lb2)& $(L,\delta)$ is a Lie coalgebra,\\
(Lb3)& For any $x,y\in L$,
$\delta([x,y])=x\cdot\delta(y)-y\cdot\delta(x)$.\end{tabular}

The compatibility condition (Lb3) shows that  $\delta$ is a
derivation map.

In the following lemmas, $c$ is an arbitrary constant.

\begin{lemm}
For any finite dimensional Lie algebra $(L,[,])$, the dual space
$L^*$ has a Lie coalgebra structure defined by
$$\delta_{L^*}(f^*)=\sum_{(f)}f_1\otimes f_2: x\otimes y\longmapsto f_1(x)f_2(y)=cf^*([x,y]).$$
\end{lemm}

\begin{lemm}
For any finite dimensional Lie coalgebra $(L,\delta)$, the dual
space $L^*$ has a Lie algebra structure defined by
$$[f^*,g^*]: x\longmapsto cf^*(x_1)g^*(x_2),$$
where $\delta(x)=\sum_{(x)}x_1\otimes x_2.$
\end{lemm}
These two lemmas are natural conclusions and easy to be verified.

\section{What is a \hxl }

\begin{defi}
Suppose that $(L,[\cdot,\cdot])$ is a Lie algebra and $(L,\delta)$
is a Lie coalgebra. A triple $(L,[\cdot,\cdot],\delta)$ is called
a \hxl if the compatibility condition
$$[\cdot,\cdot]\circ\delta=\id_L$$ holds, and $\delta$ is called a
\hx of $L$.

If in the compatibility condition, $\id_L$ is replaced by a
non-degenerate diagonal matrix, then $(L, [,], \delta)$ is called
a weak \hxl algebra and $\delta$ is called a weak \hx.
\end{defi}

\begin{rem}
Obviously, a \hxl $L$ should satisfies $[L,L]=L$.
\end{rem}

\begin{rem}
If $(L, [,], \delta)$ is  a finite dimensional (weak) \hxl algebra,
so is $(L^*, \delta^*, [,]^*)$, where \beqs
\delta^*(f\otimes g)(x)&=& (f\otimes g)\delta(x),\\
{[,]}^*(f)(x\otimes y)&=&f({[x,y]}),\eeqs for all $x,y\in L$ and
$f,g\in L^*$. This follows from the fact that $V\longrightarrow
V^*$ is a contravariant functor.
\end{rem}
\vskip5pt

\section{Co-split Lie algebras of type $A$}

Suppose that $L$ is a complex simple Lie algebra of type $A_n$, then
it can be realized as the special linear Lie algebra $sl_{n+1}(\mc)$
with basis
$$\{E_{i,j},E_{j,i},E_{i,i}-E_{j,j} \mid 1\leq i<j\leq n+1\}.$$
The Lie bracket is the commutator
$$[E_{i,j},E_{k,l}]=\delta_{j,k}E_{i,l}-\delta_{l,i}E_{k,j}.$$

Define a linear map $\delta: sl_{n+1}(\mc)\longrightarrow
sl_{n+1}(\mc)\otimes sl_{n+1}(\mc)$ as \beqs
\delta(E_{i,j})&=&\frac1{2n+2}\sum_{k=1}^{n+1}(E_{i,k}\otimes
E_{k,j}-E_{k,j}\otimes E_{i,k}). \eeqs

\begin{prop}
$\delta$ is well-defined.
\end{prop}
{\it Proof.} Assume that $i\not=j$, then\beqs
\delta(E_{i,j})&=&\frac1{2n+2}\sum_{k=1}^{n+1}(E_{i,k}\otimes
E_{k,j}-E_{k,j}\otimes E_{i,k})\\
&=&\frac1{2n+2}\sum_{k\not=i,j}(E_{i,k}\otimes
E_{k,j}-E_{k,j}\otimes E_{i,k})\\
&&+\frac1{2n+2}[(E_{i,i}-E_{j,j})\otimes E_{i,j}-E_{i,j}\otimes
(E_{i,i}-E_{j,j})]. \eeqs \beqs
\delta(E_{i,i}-E_{j,j})&=&\frac1{2n+2}\sum_{k=1}^{n+1}(E_{i,k}\otimes
E_{k,i}-E_{j,k}\otimes E_{k,j})\\
&&-\frac1{2n+2}\sum_{k=1}^{n+1}(E_{k,i}\otimes
E_{i,k}-E_{k,j}\otimes E_{j,k})\\
&=& \frac1{2n+2}\sum_{k\not=i}(E_{i,k}\otimes
E_{k,i}-E_{k,i}\otimes
E_{i,k})\\
&&-\frac1{2n+2}\sum_{k\not=j}(E_{j,k}\otimes
E_{k,j}-E_{k,j}\otimes E_{j,k}).\eeqs Hence $\delta$ is
well-defined.\hfill$\square$\vskip5pt

\begin{theo}
$(sl_{n+1}(\mc),\delta)$ is a Lie coalgebra.
\end{theo}
{\it Proof.} At first, it is clear that $(1+\tau)\circ\delta=0$. By
a direct calculation, we have \beqs
&&((1\otimes\delta)\circ\delta)(E_{i,j})\\
&=&\frac1{2n+2}(1\otimes\delta)\left(\sum_{k=1}^{n+1}(E_{i,k}\otimes
E_{k,j}-E_{k,j}\otimes E_{i,k})\right)\\
&=&\frac1{4(n+1)^2}\sum_{1\leq k,l\leq n+1}(E_{i,k}\otimes
E_{k,l}\otimes E_{l,j}-E_{i,k}\otimes E_{l,j}\otimes
E_{k,l})\\
&&-\frac1{4(n+1)^2}\sum_{1\leq k,l\leq n+1}(E_{k,j}\otimes
E_{i,l}\otimes E_{l,k}-E_{k,j}\otimes E_{l,k}\otimes E_{i,l})\\
&=&\frac1{4(n+1)^2}\sum_{1\leq k,l\leq n+1}(E_{i,k}\otimes
E_{k,l}\otimes E_{l,j}-E_{l,j}\otimes E_{i,k}\otimes
E_{k,l})\\
&&-\frac1{4(n+1)^2}\sum_{1\leq k,l\leq n+1}(E_{i,l}\otimes
E_{k,j}\otimes E_{l,k}-E_{k,j}\otimes E_{l,k}\otimes E_{i,l}). \eeqs
Hence,
$$(1+\xi+\xi^2)\circ(1\otimes\delta)\circ\delta=0,$$
that is, $\delta$ satisfies the anti-symmetriy property and the
Jacobi identity. Then $(sl_{n+1}(\mc),\delta)$ is a Lie
coalgebra.\hfill $\square$\vskip5pt

\begin{theo}
$(sl_{n+1}(\mc),[\cdot,\cdot],\delta)$ is a \hxl algebra.
\end{theo}
{\it Proof.} For $i\not=j$, it is easy to check that
\beqs
([\cdot,\cdot]\circ\delta)(E_{i,j})&=&\frac1{2n+2}\sum_{k=1}^{n+1}([E_{i,k},E_{k,j}]-[E_{k,j},E_{i,k}])\\
&=&E_{i,j},\eeqs
\beqs([\cdot,\cdot]\circ\delta)(E_{i,i}-E_{j,j})&=&\frac1{2n+2}\sum_{k=1}^{n+1}([E_{i,k},E_{k,i}]-[E_{k,i},E_{i,k}])\\
&&-\frac1{2n+2}\sum_{k=1}^{n+1}([E_{j,k},E_{k,j}]-[E_{k,j},E_{j,k}])\\
&=&E_{i,i}-E_{j,j},\eeqs that is, $[\cdot,\cdot]\circ\delta=\id$,
also by Theorem 4.1. So, the theorem holds.\hfill$\square$\vskip5pt

\section{Dual Lie algebras, Killing form and adjoint representation}

In this section, we discuss the interrelation of the Killing form
and the adjoint representation for the Lie algebra of type $A$
within our new terminology.

\begin{theo}
$((sl_n)^*,-2n\delta^*)$ is a Lie algebra isomorphic to $sl_n$,
the isomorphism is given by
$$B: f_{i,j}\longmapsto E_{j,i},$$
where $\{f_{i,j}\mid 1\leq i,j\leq n\}$ forms a basis of
$(gl_n)^*\supset(sl_n)^*$, and
$$f_{i,j}(E_{k,l})=\delta_{i,k}\delta_{j,l}.$$
\end{theo}
{\it Proof.} By definition, we have \beqs
-2n\delta^*(f_{i,j}\otimes
f_{k,l})(E_{s,t})&=&-\sum_{r=1}^n(f_{i,j}(E_{s,r})f_{k,l}(E_{r,t})-f_{i,j}(E_{r,t})f_{k,l}(E_{s,r}))\\
&=&-\delta_{j,k}(f_{i,j}(E_{s,j})f_{j,l}(E_{j,t})+\delta_{i,l}f_{i,j}(E_{i,t})f_{k,i}(E_{s,i})\\
&=&-(\delta_{j,k}f_{i,l}-\delta_{i,l}f_{k,j})(E_{s,t}),\eeqs then
$(sl_n)^*$ is a Lie algebra under bracket $-2n\delta^*$, and $B$
is an isomorphism.\hfill$\square$\vskip5pt

Define a  bilinear form $(,)_B:sl_n\times sl_n\longrightarrow \mc$
as $(x,y)_B=B^{-1}(x)(y)$.

\begin{theo}
$(,)_B$ is just a non-zero scalar of the Killing form.
\end{theo}
{\it Proof.} This result is direct.\hfill$\square$\vskip5pt Now we
can consider the following maps:
$$\xymatrix{sl_n\ar@{->}^{2n\delta\quad}[r]&sl_n\otimes sl_n\ar@{->}^{\id_{sl_n}\otimes B^{-1}\quad}[rr]&&sl_n\otimes(sl_n)^*\ar@{=>}^{\eta}[r]&\hbox{\bf End}(sl_n)}$$
where $\eta$ is an isomorphism, $\eta(x\otimes f)(y)=f(y)x$.
\begin{theo}
For the adjoint representation $$\ad:sl_n\longrightarrow \hbox{\bf
End}(sl_n),$$ we have
$$\ad=2n\eta\circ(\id_{sl_n}\otimes B^{-1})\circ\delta.$$
\end{theo}
{\it Proof.} For any $E_{i,j},E_{k,l}$, we have \beqs
2n\eta\circ(\id_{sl_n}\otimes
B^{-1})\circ\delta(E_{i,j})(E_{k,l})&=&\sum_{s=1}^n(f_{j,s}(E_{k,l})E_{i,s}-f_{s,i}(E_{k,l})E_{s,j})\\
&=&\delta_{j,k}E_{i,l}-\delta_{i,l}E_{k,j}\\
&=&\ad(E_{i,j})(E_{k,l}).\eeqs \hfill$\square$\vskip5pt

\begin{rem}
For convenience, many computations are made in $gl_n$ or
$(gl_n)^*$, but the results always hold in $sl_n$ or $(sl_n)^*$.
\end{rem}

\section{Co-splitting Theorem}

In this section, we prove the following theorem:
\begin{theo}
Any finite dimensional complex simple Lie algebra has
a \hxl structure.
\end{theo}

\noindent{\it Proof.} \ For a simple Lie algebra $L$ of type $X_l$
rather than of type $A$, our proof is divided into following steps.

\vskip.2cm {\noindent\bf Step 1:}

Suppose that $V$ is a non-trivial irreducible $X_l$-module of
dimension $n$. Then there is an injection
$$\rho: L\longrightarrow sl_n\subset\hbox{\bf End}(V),$$
and it is easy to check that the bilinear form $(,)_B$ of $sl_n$
is still non-degenerate over $\rho(L)$.

\vskip.2cm {\noindent\bf Step 2:}

Let $M$ be the orthogonal complement of $\rho(L)$ with respect to
$(,)_B$, that is,
$$M=\{m\in sl_n \mid (m,\rho(L))_B=0\}.$$
Then $M$ is a $\rho(L)$-submodule and $sl_n=\rho(L)\bigoplus M$.

\vskip.2cm {\noindent\bf Step 3:}

For any element $a=x+v\in sl_n\otimes sl_n$, if we have $x\in
\rho(L)\otimes \rho(L)$ and $v\in \rho(L)\otimes M+M\otimes
\rho(L)+M\otimes M$, the projective map $\hbox{\bf
Proj}^{sl_n\otimes sl_n}_{\rho(L)\otimes \rho(L)}$ is defined to map
$a$ to $x$. Now we write $\delta_{res}=:\hbox{\bf Proj}^{sl_n\otimes
sl_n}_{\rho(L)\otimes \rho(L)}\circ\delta|_{\rho(L)}$, where
$\delta$ is given in Section 4., then we have
\begin{lemm}
 $(\rho(L), \delta_{res})$ is a Lie coalgebra.
\end{lemm}
{\it Proof.} At first, it is easy to show that $\delta$ is an
injective map of $sl_n$-module, hence of $\rho(L)$-modules.

By Theorem 5.3, $\delta$ is equivalent to the adjoint
representation, so it is easy to know that $\delta(\rho(L))\subset
\rho(L)\otimes \rho(L)+M\otimes M\cong\rho(L)\otimes
\rho(L)^*+M\otimes M^*$. Now the skew-symmetry of $\delta_{res}$
is clear.

Furthermore, we have
$$(1\otimes\delta_{res})\circ\delta_{res}=\hbox{\bf Proj}^{sl_n\otimes
sl_n\otimes sl_n}_{\rho(L)\otimes \rho(L)\otimes
\rho(L)}\circ[(1\otimes\delta)\circ\delta]|_{\rho(L)},$$ it is
obvious by the contained relation
$$(1\otimes\delta)\circ\delta(\rho(L))\subset \rho(L)\otimes \rho(L)\otimes \rho(L)+ \rho(L)\otimes M\otimes M +M\otimes \rho(L)\otimes M+M\otimes M\otimes \rho(L),$$
thus we have proved this lemma. \hfill$\square$\vskip5pt

\vskip.2cm {\noindent\bf Step 4:}

\begin{lemm}
$$[,]\circ\delta_{res}= \hbox{\it a non-zero scalar of}\quad \id_{\rho(L)}.$$
\end{lemm}
{\it Proof.} Suppose that $\Delta^+$ is the positive root system
of $X_l$ and $\gamma$ is the highest root. It is easy to find a
basis of $\rho(L)$
$$\{X_{\pm\alpha}, h_i\mid i=1,\cdots,l; \alpha\in\Delta^+\}$$
such that $(h_i,h_j)_B=\delta_{i,j}$ and $\alpha(h_i)\in\mr$.

Since $\gamma$ is the highest root, then for any
$\alpha\in\Delta^+$, $[E_\gamma,E_\alpha]=0$. By the property of
$\delta$ (Theorem 5.3) and definition of $\delta_{res}$, we have
$$2n\delta_{res}(X_{\gamma})=\sum_{i=1}^l[X_{\gamma},h_i]\otimes
h_i+\sum_{\alpha\in\Delta^+}[X_{\gamma},X_{-\alpha}]\otimes
\frac{X_\alpha}{(X_{\alpha},X_{-\alpha})_B},$$ hence \beqs
2n{[,]}\circ\delta_{res}(X_{\gamma})&=&\sum_{i=1}^l[[X_{\gamma},h_i],h_i]
+\sum_{\alpha\in\Delta^+}\frac{[[X_{\gamma},X_{-\alpha}],X_\alpha]}{(X_{\alpha},X_{-\alpha})_B}\\
&=&\left[\sum_{i=1}^l\gamma(h_i)^2+\sum_{\alpha\in\Delta^+}\gamma(h_\alpha)\right]X_\gamma,\eeqs
where
$$h_\alpha=[X_{\alpha},X_{-\alpha}]/(X_{\alpha},X_{-\alpha})_B=\frac{2\alpha}{(\alpha,\alpha)},$$
the second assertion holds by $\rho(L)\cong\rho(L)^*$.

Clearly,
$$\sum_{i=1}^l\gamma(h_i)^2+\sum_{\alpha\in\Delta^+}\gamma(h_\alpha)>0,$$
then ${[,]}\circ\delta_{res}(X_{\gamma})\not=0$.

Secondly,  $\delta_{res}(X_\gamma)$ is a highest weight vector of
$L$-module $\rho(L)\otimes \rho(L)\cong L\otimes L$,  thus the
equation in this lemma holds.\hfill$\square$\vskip5pt

Up to now, we have completed the proof of Theorem 6.1.

We also obtain the following result.
\begin{theo}
For any type $X_l$, we have
$$2n\eta\circ(\id_{sl_n}\otimes B^{-1})\circ(\delta_{res})=\ad_{\rho(L)}|_{\rho(L)},$$
where
$\ad_{\rho(L)}|_{\rho(L)}(\rho(x)+m)=\ad(\rho(x))|_{\rho(L)}$ for
any $x\in L, m\in M$.
\end{theo}
{\it Proof.} For any $x,y\in L$, by the obvious fact
$[\rho(x),\rho(y)]\in \rho(L)$, we have \beqs
&&\left[2n\eta\circ(\id_{sl_n}\otimes
B^{-1})\circ\delta\right](\rho(x))(\rho(y))
\\&=&\left[2n\eta\circ(\id_{sl_n}\otimes
B^{-1})\circ(\delta_{res})\right](\rho(x))(\rho(y))\\
&=&\ad(\rho(x))(\rho(y))\\
&=&\ad_{\rho(L)}(\rho(x))|_{\rho(L)}(\rho(y)), \eeqs and by the
definition of $\delta_{res}$, $$[2n\eta\circ(\id_{sl_n}\otimes
B^{-1})\circ(\delta_{res})](\rho(x))(M)=0.$$ So, the claim is true.
\hfill$\square$\vskip5pt

\begin{rem} This work shows that for any finite dimensional semi-simple Lie algebra $L$
over the complex field $\mc$ (or, equivalently, over any
algebraically closed field with characteristic zero), there exists
some important relation between its Killing form and adjoint action.
Hence our new algebraic structure is proved to be very useful.
However, much more problems about it need to be solved.
\end{rem}

 \vskip10pt \centerline{\bf Acknowledgements}

The first author is supported by the Science Foundation of Jiangsu
University. He is grateful to Prof. R. Farnsteiner for his valuable
suggestion and discussion.

The second author is supported in part by the NNSF of China (Grants:
10971065, 10728102), the PCSIRT and the RFDP from the MOE of China,
the National \& Shanghai Leading Academic Discipline Projects
(Project Number: B407).

\end{document}